\documentstyle[twoside]{article}
\def\qmod#1#2{{\raise1pt\hbox{$#1$}\kern-1pt\big/
               \kern-1pt\raise-1pt\hbox{$#2$}}}

\newfam\msbfam
\font\tenmsb=msbm10 at 10 pt
\font\sevenmsb=msbm10 at 7pt
\font\fivemsb=msbm10 at 5pt

\textfont\msbfam=\tenmsb
\scriptfont\msbfam=\sevenmsb
\scriptscriptfont\msbfam=\fivemsb

\def\Bbb{\fam\msbfam\tenmsb}

\def\Z{{\Bbb Z}}

\newfam\meuffam

\font\tenmeuf=eufm10
\font\sevenmeuf=eufm7
\font\fivemeuf=eufm5

\textfont\meuffam=\tenmeuf
\scriptfont\meuffam=\sevenmeuf
\scriptscriptfont\meuffam=\fivemeuf

\def\germ{\fam\meuffam\tenmeuf}

\def\g{{\germ g}}

\def\Pg{{\germ P}}

\def\qed {\hfill\vrule height6pt width6pt depth0pt \smallskip}

\def\ra{\rightarrow}
\def\textmap#1{\mathop{\vbox{\ialign{
                                ##\crcr
    ${\scriptstyle\hfil\;\;#1\;\;\hfil}$\crcr
    \noalign{\kern-1pt\nointerlineskip}
    \rightarrowfill\crcr}}\;}}

\def\textlmap#1{\mathop{\vbox{\ialign{
                                ##\crcr
    ${\scriptstyle\hfil\;\;#1\;\;\hfil}$\crcr
    \noalign{\kern-1pt\nointerlineskip}
    \leftarrowfill\crcr}}\;}}

\newfam\meuffam
\font\tenmeuf=eufm10
\font\sevenmeuf=eufm7
\font\fivemeuf=eufm5

\textfont\meuffam=\tenmeuf
\scriptfont\meuffam=\sevenmeuf
\scriptscriptfont\meuffam=\fivemeuf

\def\germ{\fam\meuffam\tenmeuf}

\def\g{{\germ g}}

\def\picture#1by#2(#3){
\vbox to #2 {
  \hrule width #1 height 0pt depth 0pt \vfill \special{picture #3}}
}

\def\scaledpicture#1by#2(#3scaled#4){{
\dimen0=#1  \dimen1=#2
\divide\dimen0 by 1000 \multiply\dimen0 by #4
\divide\dimen1 by 1000 \multiply\dimen1 by #4
\picture \dimen0 by \dimen1 (#3 scaled #4)}}
\def\dfigure#1by#2(#3scaled#4offset#5:#6)
  {\medskip
   \vglue 2mm minus 2mm
   $$
     \hbox{
       \hglue#5
       {\scaledpicture #1 by #2 (#3 scaled #4)}
     }
   $$
   \par\goodbreak
   \vglue 2mm minus 2mm
   \medskip}

\begin{document}

\def\Pr{{\rm Pr}}
\def\tr{{\rm Tr}}
\def\coker{{\rm coker}}
\def\ad{{\rm ad}}
\def\End{{\rm End}}
\def\Pic{{\rm Pic}}
\def\NS{{\rm NS}}
\def\deg{{\rm deg}}
\def\Hom{{\rm Hom}}
\def\Aut{{\rm Aut}}
\def\Herm{{\rm Herm}}
\def\Vol{{\rm Vol}}
\def\pf{{\bf Proof: }}
\def\id{{\rm id}}
\def\im{{\rm im}}
\def\rk{{\rm rk}}
\def\Sp{{\rm Sp}}
\def\Spin{{\rm Spin}}
\def\h{{\bf H}}
\def\dv{\bar\partial}
\def\dva{\bar\partial_A}
\def\da{\partial_A}
\def\p{\partial\bar\partial}
\def\pa{\partial_A\bar\partial_A}
\def\sw{Seiberg-Witten}
\def\proj{{\rm pr}}
\def\ub{\underbar}
\def\gr{{\scriptscriptstyle|}\hskip -4pt{\g}}
\def\Dr{{\raisebox{0.15ex}{$\not$}}{\hskip -1pt {D}}}
\def\dr{{\raisebox{0.15ex}{${\scriptstyle\not}$}}{\hskip -0.8pt
{\scriptstyle D}}}
\def\subsetint{{\  {\subset}\hskip -2.45mm{\raisebox{.28ex}
{$\scriptscriptstyle\subset$}}\ }}
\def\nr{\parallel}
\newtheorem{sz}{Satz}
\newtheorem{szfr}{Satzfr}
\newtheorem{thry}[sz]{Theorem}
\newtheorem{thfr}[szfr]{Th\'eor\`eme}
\newtheorem{pr}[sz]{Proposition}
\newtheorem{re}[sz]{Remark}
\newtheorem{co}[sz]{Corollary}
\newtheorem{cofr}[szfr]{Corollaire}
\newtheorem{dt}[sz]{Definition}
\newtheorem{lm}[sz]{Lemma}
\noindent G\'eom\'etrie Analytique /{\sl Analytical Geometry}
\\ \\
\centerline{\Large{\bf Holomorphic vector bundles on non-algebraic surfaces}}\vspace{2mm}
\vspace*{2mm}
\centerline{\large   Andrei
Teleman \hspace{3cm} Matei Toma }\\  
{\footnotesize{\sl Abstract --}} {\footnotesize    The existence problem for holomorphic
structures on vector bundles over non-algebraic surfaces is in general still open. We
solve this problem in the case of rank 2 vector bundles over   K3 surfaces and in the
case of vector bundles of arbitrary rank over all  known surfaces of class VII. Our
methods, which are based on Donaldson theory and deformation theory, can be used to
solve the existence problem of holomorphic vector bundles on further classes of
non-algebraic surfaces.}
\vskip 6pt

\centerline{{\bf Fibr\'es vectoriels holomorphes sur les surfaces non-alg\'ebriques}}
\vskip 4pt
{\footnotesize{\sl R\'esum\'e --}} {\footnotesize  Le probl\`eme de l'existence des structures
holomorphes sur les fibr\'es vectoriels au-dessus des surfaces non-alg\'ebriques est en g\'en\'eral
encore ouvert. Nous r\'esolvons ce probl\`eme pour les fibr\'es de rang 2 sur les surfaces K3 et pour les
fibr\'es   de rang arbitraire sur toutes les surfaces connues de la classe VII. Nos m\'ethodes, qui
s'appuient sur la th\'eorie de Donaldson et sur la th\'eorie des d\'eformations, peuvent
\^etre utilis\'ees pour r\'esoudre le probl\`eme de l'existence des fibr\'es vectoriels holomorphes sur 
d'autres classes de surfaces non-alg\'ebriques.}
\vskip 8pt

{\sl Version fran\c caise abr\'eg\'ee } -   Soit $E$ un fibr\'e vectoriel complexe sur
une surface complexe compacte
$X$. Nous posons
$$\Delta(E):=2\rk(E) c_2(E)- (\rk(E)-1) c_1(E)^2\ .$$

Pour tout $a\in \NS(X)$ et pour tout entier positif $r$ soit
$$m(r,a):=r\inf \left\lbrace-\sum_{i=1}^{r}\left(\frac{a}{r}-
\mu_i\right)^2\left|\right.\mu_1,\ldots,\mu_r\in\NS(X)\mbox{\ avec\ }
\sum_{i=1}^{r}\mu_i=a\right\rbrace .$$
Ce nombre est un entier non-negatif si $X$ est non-alg\'ebrique 
[2].  Un fibr\'e holomorphe ${\cal E}$ s'appelle {\it
filtrable} s'il existe une filtration   $0\subset {\cal
E}_1\subset \dots\subset {\cal E}_{r-1}\subset {\cal E}_r={\cal E}$ avec des faisceaux
coherents ${\cal E}_i$ tels que $\rk({\cal E}_i)=i$  pour tous les $i$. Le
probl\`eme de l'existence des structures holomorphes filtrables sur les  
surfaces complexes a  \'et\'e compl\'etement r\'esolu dans [2]. Des r\'esultats sur les
structures non-filtrables ont  \'et\'e obtenus dans   [1], [12-14]. Notre premier
r\'esultat est
\vskip 5pt
{\sc Th\'eor\`eme 1.} {\it Soit $E$ un fibr\'e complexe de rang 2 diff\'erentiable sur
une surface    K3  non-alg\'ebrique $X$.  Alors, \`a l'exception du cas
$a(X)=0$, $\Delta(E)=4$ et
$c_1(E)\in 2 NS(X)$, on a \\
i) $E$ admet une structure holomorphe filtrable   si et seulement si
$$c_1(E)\in NS(X)\ {\rm et}\ \Delta(E)\geq  m(2,c_1(E))\ .$$
ii) $E$ admet une structure holomorphe  si et seulement si
$$c_1(E)\in NS(X)\ {\rm et}\ \Delta(E)\geq \min(6,m(2,c_1(E)))\ ,$$
Dans le cas except\'e, $E$ n'admet aucune structure  holomorphe.
}\vspace{1.5mm}\\
Pour d\'emontrer l'existence des structures holomorphes nous utilisons la
non-trivialit\'e de  l'invariant  polynomial de Donaldson $q_P(X)$
associ\'e   au 
$PU(2)$-fibr\'e $P$ correspondant \`a une structure hermitienne sur $E$.
\vskip 5pt
{\sc Remarque:} La th\'eorie de Donaldson permet de r\'esoudre le   probl\`eme
de l'existence des structures holomorphes dans le cas
$\rk(E)=2$   pour les surfaces k\"ahleriennes  $X$ telles que
$H^0(K_X)\textmap{\otimes^2} H^0(K_X^{\otimes2})$ soit surjective.  Pour une telle
surface $X$, fixons
$c\in NS(X)$ et d\'esignons par $\bar c\in H^2(X,\Z_2)$ sa reduction mod
2. On peut montrer que 
$$\min\{\Delta(E)|\ c_1(E)=c,\ E\ {\rm admet\ des\ structures\
holomorphes}\}=$$
$$\min\{-p_1(P)|\ w_2(P)=\bar c,\ q_P(X)\ne 0 \}\ ,
$$
lorsque la valeur minimale \`a droite
est inf\'erieure \`a $m(c_1,2)$.

Ce r\'esultat montre que, pour cette  classe de surfaces k\"ahleriennes,
l'apparition des structures holomorphes dans le cas
$\Delta(E)<m(c_1(E),2)$ est d\'etermin\'ee par {\it le type topologique
diff\'erentiable} de la 4-vari\'et\'e diff\'erentiable sous-jacente.
\vskip 5pt
Nous dirons qu'une   surface complexe $X$ a la propri\'et\'e $\Pg_r$ si  tout fibr\'e
vectoriel de rang $r$ sur $X$ qui admet une  structure holomorphe,    admet aussi une
structure  filtrable.   
\vskip 3pt
{\sc Proposition. } {\it \label{inegalitati}
Soit $X$ une surface complexe , $\pi :\tilde X \rightarrow X$ une
modification propre et $r \geq 2$ un  entier. Alors   $\tilde X$  a la
propri\'et\'e
$\Pg_r$ si et seulement si $X$ a cette propri\'et\'e.  }
\vspace{1.5mm}\\
L'id\'ee de la  d\'emonstration est de montrer l'in\'egalit\'e
$$
\Delta({\cal E})-m(r,c_1({\cal E}))\geq\Delta(\pi_*{\cal E})-
m(r,c_1(\pi_*{\cal E}))\ .
$$
pour tout fibr\'e holomorphe ${\cal E}$ de rang $r$ sur $\tilde X$.
\vskip 5pt
Notre deuxi\`eme r\'esultat r\'esout le probl\`eme de l'existence des structures
holomorphes sur les fibr\'es vectoriels au-dessus de  {\it toutes} les surfaces
{\it connues} de la classe $VII$:
\vskip 5pt
{\sc Th\'eor\`eme 2.} {\it Soit $X$ une surface de la classe VII  dont le   mod\`ele 
minimal est soit une surface \`a
$b_2=0$ soit  une surface  avec un cycle de courbes
rationnelles. Soit $E$ un fibr\'e vectoriel diff\'erentiable
de rang $r$ sur    $X$.  Les propri\'et\'es suivantes sont \'equivalentes:\\
1.  $E$ admet une  structure holomorphe, \\  
2. $E$ admet une  structure holomorphe filtrable,\\
3. $c_1(E)\in\NS(X)\mbox{\ et}
     \    \Delta(E)\geq m(r,c_1(E)) .$}%
\vspace{1.5mm}\\
Lorsque le   mod\`ele 
minimal est une surface \`a $b_2=0$, on applique la Proposition ci-dessus et le
r\'esultat fondamental de [2]. Dans le deuxi\`eme cas l'id\'ee de la d\'emonstration est
de r\'eduire le probl\`eme par petite d\'eformation
au cas d'une surface de Hopf \'eclat\'ee, et d'appliquer de nouveau la m\^eme
Proposition.
\\ \noindent\hbox to 5cm{\hrulefill} \\

1. {\sc Notations}. -  Let $E$ be a differentiable complex vector bundle over a complex surface
$X$. The discriminant of $E$ is defined by
$$\Delta(E):=2\rk(E) c_2(E)- (\rk(E)-1) c_1(E)^2\ .$$
For each $a\in \NS(X)$ and any positive integer $r$ we put
$$m(r,a):=r\inf \left\lbrace-\sum_{i=1}^{r}\left(\frac{a}{r}-
\mu_i\right)^2\left|\right.\mu_1,\ldots,\mu_r\in\NS(X)\mbox{\ with\ }
\sum_{i=1}^{r}\mu_i=a\right\rbrace .$$
Note that $m(r,a)$ equals $-\infty$ if $X$ is algebraic and is a non-negative integer if
$X$ is non-algebraic [2].  We recall that a holomorphic vector bundle ${\cal E}$ is called
{\it filtrable} if there exists a filtration $0\subset {\cal E}_1\subset \dots \subset
{\cal E}_{r-1}\subset {\cal E}_r={\cal E}$ with coherent sheaves ${\cal E}_i$ such
that $\rk({\cal E}_i)=i$ for all $i$. The existence problem for filtrable
holomorphic  structures on complex surfaces was completely solved in [2]. For
non-filtrable structures the existence problem was solved over tori for rank 2 and over
primary Kodaira surfaces for arbitrary rank (see [12-14], [1]). We treat here the case of
K3 surfaces and of class VII surfaces.
\vskip 5pt

2.   {\sc The case of K3 surfaces}. -    The following result answers completely the existence  question
for holomorphic structures on rank 2 vector bundles over K3 surfaces.
\vskip 5pt 
{\sc Theorem  1.} {\it Let $E$ be a rank 2 differentiable  vector bundle over a
non-algebraic K3 surface $X$.  Excepting the case when $a(X)=0$, $\Delta(E)=4$ and
$c_1(E)\in 2 NS(X)$, one has \\
i) $E$ admits a filtrable holomorphic structure if and only if 
$$c_1(E)\in NS(X)\ {\rm and}\ \Delta(E)\geq  m(2,c_1(E))\ .$$
ii) $E$ admits a holomorphic structure if and only if 
$$c_1(E)\in NS(X)\ {\rm and}\ \Delta(E)\geq \min(6,m(2,c_1(E)))\ ,$$
In the excepted case, $E$ admits no holomorphic structure.
}
\vspace{1.5mm}\\
\pf  The first statement is a particular case of the main theorem in [2], where
also the non-existence of holomorphic structures in the excepted case is proved.  It
remains to show:\\ A. If $\Delta(E)< \min(6,m(2,c_1(E)))$, then $E$ admits no
holomorphic structure,\\ B.  If  $c_1(E)\in NS(X)$, and $6\leq\Delta(E)<
m(2,c_1(E))$, then $E$ does admit holomorphic structures.

To prove  A, notice first that, if $\Delta(E)<  m(2,c_1(E))$, then any holomorphic
structure   on $E$ is non-filtrable (by $i))$, hence simple.  Let ${\cal E}$ be such
such a holomorphic structure.  Since
$K_X={\cal O}_X$, one gets
$H^0({\cal E}nd_0({\cal E}))=H^2({\cal E}nd_0({\cal E}))=0$, and, by Riemann-Roch,
this yields $h^1({\cal E}nd_0({\cal E}))<0$ when $\Delta(E)<6$.

To prove B, let $\delta$ be a fixed a holomorphic connection on $\det(E)$ and fix a hermitian
metric $h$ on $E$.    
Denote also by $P$ the $PU(2)$-bundle associated with the $U(2)$-bundle $(E,h)$. Its
Pontrjagin class is $p_1(P)=-\Delta(E)$. The Kobayashi-Hitchin correspondence ([3], [7]) 
gives an isomorphism of  real analytic spaces
$${\cal M}^{\rm ASD}_g(P)\simeq {\cal M}^{\rm pst}_{g,\delta}(E)\ ,
$$
for every Gauduchon metric $g$ on $X$. Here ${\cal M}^{\rm ASD}_g(P)$ denotes the moduli space
of  $g$ - ASD connections on $P$ and ${\cal M}^{\rm st}_{g,\delta}(E)$
denotes the moduli space of
$g$ - poly-stable holomorphic structures on $E$ which induce $\delta$ on $\det(E)$.  In our
case we have
$w_2(P)\ne 0$. Indeed, if $w_2(P)$ vanished, then $c_1(E)\in 2 H^2(X,\Z)\cap
NS(X)=2NS(X)$, which would imply $m(2,c_1(E))=0$.  

Since $w_2(P)\ne 0$ and $X$ is simply connected, the Donaldson polynomial
invariant associated with
$P$ is well-defined [3]\footnote{The Donaldson polynomial invariants
associated with
$PU(2)$-bundles with $w_2=0$ on simply connected manifolds are defined only for sufficiently
large instanton number (the stable range).  A more refined theory [4], which uses the
thickened ASD moduli space leads to well defined invariants for any value of $c_2$, but this
theory is not useful for our purposes.}.  On  the other hand, by the same argument  as above,
 $h^0({\cal E}nd_0({\cal E}))=h^2({\cal E}nd_0({\cal E}))=0$  for any bundle ${\cal
E}$ with $\Delta({\cal E})<
m(2,c_1({\cal E}))$, hence the moduli space
${\cal M}^{\rm ASD}_g(P)$ and all similar moduli spaces corresponding to the lower Uhlenbeck
strata, are   regular of expected dimension.  Therefore, one can compute the Donaldson
polynomial invariant corresponding to $P$ using the Uhlenbeck compactification
$\overline{{\cal M}}^{\rm ASD}_g(P)$.  But, by classical results in gauge theory
([3], [6],[9]) the Donaldson polynomial invariant $q_P(X)$  associated to any 
$PU(2)$ bundle 
$P$ with
$-p_1(P)\geq 6$  is non-trivial. This shows that ${\cal M}^{\rm ASD}_g(P)$ (hence also
${\cal M}^{\rm st}_{g,\delta}(E)$) cannot be empty.
\qed
\vskip 3pt
{\sc Remark:} Donaldson theory can be used to solve the existence problem for
$\rk(E)=2$ and all K\"ahlerian surfaces $X$ with $H^0(K_X)\textmap{\otimes^2}
H^0(K_X^{\otimes2})$ surjective. This condition ensures that any non-filtrable
holomorphic bundle defines a smooth point in the corresponding moduli space. 
For such a surface
$X$, let us fix
$c\in NS(X)$ and denote by $\bar c\in H^2(X,\Z_2)$   its reduction   mod
2. One can show that 
$$\min\{\Delta(E)|\ c_1(E)=c,\ E\ {\rm admits\ holomorphic\ structures\
}\}=$$
$$\min\{-p_1(P)|\ w_2(P)=\bar c,\ q_P(X)\ne 0 \}\ ,
$$
when the minimal value on the right is smaller than $m(c_1,2)$.  

This statement
shows that, for this class of surfaces, the existence of holomorphic
structures in the range $\Delta(E)<m(c_1(E),2)$ can be decided in terms of
the {\it differential topological type} of the underlying
differentiable 4-manifold.  The proof is based on Donaldson's
non-vanishing result for K\"ahlerian surfaces ([3], p. 378).
\vskip 5pt
3.   {\sc Blow up inequalities}. -  We say that a complex surface $X$ has the property
$\Pg_r$ if   every rank $r$  vector bundle on $X$ which admits a holomorphic structure
also admits a  filtrable one.  

By [2] a complex surface   $X$ has the property $\Pg_r$ if and only if any torsion
free sheaf ${\cal F}$ on $X$ has $\Delta({\cal F})\geq m(r,c_1({\cal F}))$.
\vskip 3pt
{\sc Proposition. } {\it \label{inegalitati}
Let $X$ be a compact complex surface, $\pi :\tilde X \rightarrow X$ a proper
modification and $r \geq 2$ an integer. Then   $\tilde X$  has the property $\Pg_r$
if and only if $X$ does.  }
\vspace{1.5mm}\\
\pf Let $X$ be non-algebraic, otherwise the statement
is trivial. We may suppose that
$\pi$ consists of a single blow-up, and let $D$
be the exceptional divisor.
Suppose   that $\Pg_r$ holds for $X$.  
By [2] any  rank $r$ holomorphic   bundle ${\cal F}$
on $X$ has $\Delta({\cal F})\geq m(r,c_1({\cal F})) $.

We will show that every holomorphic vector bundle ${\cal E}$ of rank $r$ over
$\tilde X$ satisfies $\Delta({\cal E})\geq m(r,c_1({\cal E}))$, which will
imply the property $\Pg_r$ for $\tilde X$, again by [2].

We have 
$$c_1({\cal E})=\pi^*c_1(\pi_*({\cal E}))+k [D]\ ,\eqno{(1)}$$
 where $k:=-[D]\cdot  c_1({\cal
E}) \in
\Z$.  By tensorizing ${\cal E}$ by ${\cal O}(lD), l
\in \Z$ if necessary, we may suppose that $0\leq k<r$.   The
Riemann-Roch formula   yields:
$$ \chi({\cal E}) = r\left[\chi({\cal O}_{\tilde X}) + \frac {c_1({\cal E})\cdot
c_1(\tilde X)} {2r} +\frac {c_1({\cal E})^2-\Delta({\cal E})}{2r^2} \right].$$
Since $\chi ({\cal E})=\chi(\pi_*{\cal E}) -\chi(R^1\pi_*{\cal E})$ we have
$\chi ({\cal E})\leq\chi(\pi_*{\cal E})$, which, by (1), is equivalent to
$$ \Delta({\cal E})\geq \Delta(\pi_*{\cal E})+ {k(r-k)}\ .\eqno{(2)}$$
On the other hand  
$$
m(r,c_1(\pi_*{\cal E}))= -r\max\left\lbrace\sum_{i=1}^{r}
\left(\frac{c_1(\pi_*{\cal E})}{r}-
\mu_i\right)^2\left|\right.\mu_1,\ldots,\mu_r\in\NS(X)\mbox{\ ,\ }
\sum_{i=1}^{r}\mu_i=c_1(\pi_*{\cal E})\right\rbrace $$
$$
= -r\sum_{i=1}^{r}
\left(\frac{c_1(\pi_*{\cal E})}{r}-
\mu_i^0\right)^2 \ ,$$
where the last equality holds for those $\mu_i^0$ for which  the maximum  is attained.

But
$$
  m(r,c_1({\cal E}))=-r\max\left\lbrace\sum_{i=1}^{r}
\left(\frac{c_1({\cal E})}{r}-
\nu_i\right)^2\left|\right.\nu_1,\ldots,\nu_r\in\NS(\tilde X)\mbox{\ with\ }
\sum_{i=1}^{r}\nu_i=c_1({\cal E})\right\rbrace\leq
$$
$$
\leq -r\sum_{i=1}^{k}\bigl(\frac{\pi^*c_1(\pi_*{\cal
E})+k[D]}{r}-\pi^*\mu_i^0-[D]
\bigr)^2
 -r\sum_{i=k+1}^{r}\bigl(\frac{\pi^*c_1(\pi_*{\cal E})+k[D]}{r}
-\pi^*\mu_i^0\bigr)^2= 
$$
$$
=-r\sum_{i=1}^{r}\bigl(\frac{c_1(\pi_*{\cal E})}{r}-\mu_i^0\bigr)^2+
 {k(r-k)}\ .
$$
Combining this  inequality with (2) we obtain
$$
\Delta({\cal E})-m(r,c_1({\cal E})\geq\Delta(\pi_*{\cal E})- m(r,c_1(\pi_*{\cal E}))\ .
$$
Therefore, since $\Pg_r$ holds for $X$,
$$\Delta({\cal E})\geq m(r,c_1({\cal E}))\ , $$
which proves one implication of the  Proposition.

The converse is easier;
if for a locally free sheaf ${\cal F}$ of rank $r$ on $X$ we had
$$
\Delta({\cal F})< m(r,c_1({\cal F})),
$$
then
$$
\Delta(\pi^*{\cal F})=\Delta({\cal F})< m(r,c_1({\cal F}))=m(r,c_1(\pi^*{\cal F}))
$$
proving our claim. \qed 

4. {\sc The case of class VII surfaces}. -  A compact complex surface is said to be
in {\it class $VII$} if its first  Betti number $b_1$ equals one. The subclass of such
surfaces with
$b_2=0$ is completely understood; cf. [5], [11]. 
Further surfaces in class $VII$ can be obtained by an explicit construction   due
to M. Kato. 
All known minimal surfaces in  class $VII$ with $b_2>0$ 
contain one or two cycles of rational curves. On the other hand it was shown
by Nakamura, [8], that any surface in  class $VII$ containing a 
cycle
of rational curves is a deformation of a blown-up Hopf surface. 
Using this result we are able to prove
\vskip 5pt
{\sc Theorem 2} {\it Let $X$ be a class VII surface whose minimal model  either has
$b_2=0$ or  contains a cycle of rational
curves, and let $E$ be a rank
$r$ differentiable vector bundle over   $X$.  Then the following are equivalent\\
1.  $E$ admits a holomorphic structure, \\  
2. $E$  admits a filtrable holomorphic structure,\\
3. $c_1(E)\in\NS(X)\mbox{\ and}
     \    \Delta(E)\geq m(r,c_1(E)) .$}%
\vspace{1.5mm}\\
\pf
We prove only  $1.\Rightarrow 2.$. The rest follows from [2]. By the previous
Proposition, we may suppose that
$X$ is minimal. If
$b_2(X)=0$, the statement follows again from [2]. In particular it holds for Hopf
surfaces, hence also for blown-up Hopf surfaces.

Suppose that $X$ contains a cycle of rational curves.   Let ${\cal
E}$ be a non-filtrable
 holomorphic structure on $E$. We may suppose that ${\cal E}$ does not admit any
non-trivial coherent subsheaf ${\cal F}$ of lower rank (if not, we replace ${\cal E}$ by
${\cal F}^{\vee\vee}$ and we argue by induction). 

Denote by $\Theta_X$ the tangent sheaf of  $X$.
Consider the sheaf $\Sigma$ arising as middle term of the Atiyah sequence:
$$0\ra {\cal E} nd({\cal E})\ra \Sigma \ra \Theta_X \ra 0.$$
For the deformation theory of the  pair $(X,{\cal E})$, the spaces $H^0(X,\Sigma)$,
$H^1(X,\Sigma)$, 
$H^2(X,\Sigma)$ play the roles of tangent space to the
group of automorphisms,  tangent space to the versal deformation germ
and space of obstructions to deformation respectively  (see  [10]). Moreover,
the natural map
$H^1(X,\Sigma)  \ra H^1(X,\Theta_X)$
is the tangent map at $(X,{\cal E} )$ of the natural morphism between the versal
deformation spaces of $(X,{\cal E})$ and of $X$. 

Note that $H^2(X, {\cal E} nd({\cal E}))=H^0(X, {\cal E} nd({\cal E}) \otimes
K_X)^*=0$.  Indeed, if $f\in H^0(X, {\cal E} nd({\cal E}))$, then $\det(f)\in
H^0(X,K_X^{\otimes r})$ vanishes because ${\rm kod}(X)=-\infty$, hence $\ker f$ would
contradict our assumption on ${\cal E}$. 

We also have $H^2(X,\Theta_X)=0$ (see [8]), hence
$H^2(X,\Sigma)=0$. Therefore the versal deformation space is smooth at
$(X,{\cal E})$ and the germ of the map to the versal deformation of $X$ is submersive
at $(X,{\cal E})$. In particular, over  any   deformation $X_s$ of
$X$ sufficiently close to $X$, there exists a holomorphic vector bundle ${\cal E}_s$
which is a deformation of
${\cal E}$. Choosing
$X_s$ to be a blown-up Hopf surface and 
remembering that $NS(X)=H^2(X, \Z)=H^2(X_s, \Z)=NS(X_s)$ we see that $E_s$
and hence $E$ satisfies the third condition of our theorem, hence   also the second by
[2].
\qed

\parindent0cm
\vspace{0.3cm}  
{{\bf References}}\vskip 7pt
{\small

[1]   Aprodu, M., Br\^\i nz\u anescu, V.,  Toma, M.:
{\it Holomorphic  vector bundles on primary Kodaira
surfaces}, to appear in Math. Z.

[2]   B\u{a}nic\u{a}, C., Le Potier, X.: {\it Sur l'existence
des fibr\'es vectoriels holomorphes sur les surfaces non-alg\'ebriques},
J. Reine Angew. Math. 378 (1987), 1-31.
 
[3]   Donaldson, S.,   Kronheimer, P.: {\it The Geometry of
Four-Manifolds}, Oxford Univ. Press,  1990.

[4]  Friedman, R., Morgan, J.: {\it Smooth 4-manifolds and complex surfaces},
Springer-Verlag, 1994

 [5]  Inoue, M.: {\it On surfaces of class VII$_0$},
Invent. Math. 24 (1974), 269-310.  

[6] Kronheimer, P., Mrowka, T.: {\it Embedded surfaces and the structure of Donaldson's
polynomial invariants}, J. Differential Geom. 41 (1995),  573-734.

[7]  L\"ubke, M., Teleman, A.: {\it The Kobayashi-Hitchin 
correspondence},    World Sci. Publ.,  1995.

[8]  Nakamura, I.: {\it On surfaces of class $VII_0$ with curves, II},
T\^ohoku Math. J. 42 (1990), 475-516.

[9]  O'Grady, K. G.:, I.: {\it Donaldson polynomials for K3 surfaces},
J. Diff. Geom. 35 (1992), 415-427.

[10]  Schumacher, G.: Toma, M.: 
{\it  Moduli of Kaehler manifolds equipped with Hermite-Einstein
vector bundles}, Rev. Roum. Math. Pures Appl. 38 (1993), 703-719.

[11]  Teleman, A. D.: {\it Projectively flat surfaces and 
Bogomolov's theorem on class $VII\sb 0$ surfaces},
Int. J. Math. 5 (1994), 253-264. 

[12]  Toma, M.: {\it Une classe de fibr\'es vectoriels holomorphes 
sur les
2-tores complexes}, C.R. Acad. Sci. Paris, 311 (1990), 257-258. 
 
[13]  Toma, M.: {\it Stable bundles on non-algebraic surfaces  giving 
 rise to
compact moduli spaces}, C.R. Acad. Sci. Paris, 323 (1996), 501-505.

[14]  Toma, M.: {\it Stable bundles with small $c_2$ over 2-dimensional
complex tori}, Math. Z. 232 (1999), 511-525.
\vspace{0.2cm}\\

 A. Teleman: CMI, Universit\'e de Provence,
39, Rue F. Joliot Curie, 13453 Marseille Cedex 13, France 
 and Faculty of Mathematics, University of
 Bucharest,  e-mail: teleman@cmi.univ-mrs.fr, Tel:
-33-(0)491113606, Fax: -33-(0)491-113552
\vspace{0.2cm}\\
M. Toma: Fachbereich Mathematik-Informatik, Universit\"at Osnabr\"uck,
49069 Osnabr\"uck, Germany and Mathematical Institute of the Romanian
Academy; e-mail: matei@mathematik.uni-osnabrueck.de Tel:
-49-541-9692567. Fax: -49-541-9692770. 
   
 }

\end{document}